# Двойственные подходы к задачам минимизации сильно выпуклых функционалов простой структуры при аффинных ограничениях

*Аникин А.С. (ИДСТУ СО РАН),*
*Гасников А.В. (ИППИ РАН; ПреМоЛаб ФУПМ МФТИ),*
*Двуреченский П.Е. (WIAS Berlin, ИППИ РАН),*
*Тюрин А.И. (НИУ ВШЭ),*
*Чернов А.В. (ФУПМ МФТИ)*

**Аннотация**
В статье рассматривается задача минимизации сильно выпуклой функции простой структуры (например, сепарабельной) при аффинных ограничениях. Строится двойственная задача. Для решения двойственной задачи предлагается использовать быстрый градиентный метод. В работе устанавливаются необходимые свойства этого метода, которые позволяют при весьма общих условиях восстанавливать по генерируемой этим методом последовательности в двойственном пространстве решение прямой задачи с той же точностью, что и двойственной. Несмотря на кажущуюся естественность такого подхода, стоит заметить, что в данной работе приведено решение ряда ранее неопубликованных и местами довольно тонких моментов, необходимых для строгого и полного теоретического обоснования отмеченного подхода в нужной общности.
**Ключевые слова:** прямо-двойственные методы, быстрый градиентный метод, двойственная задача, регуляризация двойственной задачи, техника рестартов, сильная выпуклость, задача PageRank.

## 1. Введение

В данной работе мы существенно обобщаем результаты недавней работы [1] (в том числе улучшаем оценку скорости сходимости метода), в которой был предложен способ решения задачи энтропийно-линейного программирования (ЭЛП) с помощью решения специальным образом регуляризованной двойственной задачи, и восстановлении по явным формулам решения прямой задачи, исходя из полученного приближенного решения двойственной задачи. Настоящая работа развивает идеи цикла статей А.С. Немировского, Ю.Е. Нестерова с соавторами и их последователей [2–13], в которых предлагались различные прямо-двойственные методы для широкого класса задач. Название "прямо-двойственные методы" было предложено в диссертации [5] для методов, которые позволяют, исходя из решения прямой (двойственной) задачи, без существенных дополнительных затрат восстанавливать (с той же точностью) решение соответствующей двойственной (прямой) задачи. Для удобства чтения статьи мы постарались привести в статье все необходимые выкладки, хотя в ряде случаев они и не являются оригинальными.

В п. 2, следуя работе [14], мы приводим быстрый градиент метод. В отличие от [1] мы также исследуем его прямо-двойственность [2] и следующее свойство: генерируемая методом последовательность точек лежит в шаре с центром в решении и радиуса равного расстоянию от точки старта метода до решения задачи. Оба эти свойства необходимы в п. 3, чтобы обосновать способ восстановления решения прямой задачи (минимизации сильно выпуклой функции простой структуры при аффинных ограничениях), исходя из сгенерированной методом последовательности в двойственном пространстве. В конце пункта мы описываем непосредственное обобщение конструкции работы [1], связанной с регуляризацией двойственной задачи. Это обобщение не требует прямо-двойственности метода, но приводит к некоторым потерям в оценках скорости сходимости. В частности, предложенный в данной статье метод, примененный к задаче ЭЛП, находит на несколько порядков быстрее решения различных задач ЭЛП, по сравнению с методом из [1].

## 2. Прямо-двойственность быстрого градиентного метода

Рассмотрим задачу выпуклой оптимизации
$$f(x) \to \min_x. \qquad (1)$$
Под решением этой задачи будем понимать такой $\bar{x}^N$, что
$$f(\bar{x}^N) - f_* \leq \varepsilon,$$
где $f_* = f(x_*)$ – оптимальное значение функционала в задаче (1), $x_*$ – решение задачи (1). Определим множество
$$B_R(x_*) = \{x : \|x - x_*\|_2 \leq R\}.$$
Пусть
$$x^{k+1} = x^k - h \nabla f(x^k), \qquad (2)$$
и при $x \in B_{\sqrt{2}R}(x_*)$, где
$$R = \|x^0 - x_*\|_2 = \|x_*\|_2,$$
выполняется условие
$$\|\nabla f(x)\|_2 \leq M.$$
Тогда из (2) с учетом этого неравенства и (4) (см. ниже), имеем
$$\|x - x^{k+1}\|_2^2 = \|x - x^k + h \nabla f(x^k)\|_2^2 = \|x - x^k\|_2^2 + 2h \langle \nabla f(x^k), x - x^k \rangle + h^2 \|\nabla f(x^k)\|_2^2 \leq$$
$$\leq \|x - x^k\|_2^2 + 2h \langle \nabla f(x^k), x - x^k \rangle + h^2 M^2.$$
Отсюда (при $x = x_*$) следует, что

$$f\left(\frac{1}{N}\sum_{k=0}^{N-1} x^k\right) - f_* \leq \frac{1}{N}\sum_{k=0}^{N-1} f(x^k) - f(x_*) \leq \frac{1}{N}\sum_{k=0}^{N-1} \langle \nabla f(x^k), x^k - x_* \rangle \leq$$

$$\leq \frac{1}{2hN}\sum_{k=0}^{N-1}\left\{\|x_* - x^k\|_2^2 - \|x_* - x^{k+1}\|_2^2\right\} + \frac{hM^2}{2} = \frac{1}{2hN}\left(\|x_* - x^0\|_2^2 - \|x_* - x^N\|_2^2\right) + \frac{hM^2}{2}.$$

Выбирая
$$h = \frac{R}{M\sqrt{N}}$$
и полагая
$$\bar{x}^N = \frac{1}{N}\sum_{k=0}^{N-1} x^k,$$
получим
$$f(\bar{x}^N) - f_* \leq \frac{MR}{\sqrt{N}}. \qquad (3)$$
Заметим, что

$$0 \le \frac{1}{2hk}\left(\|x_* - x^0\|_2^2 - \|x_* - x^k\|_2^2\right) + \frac{hM^2}{2},$$

Поэтому при $k = 0,...,N$

$$\|x_* - x^k\|_2^2 \le \|x_* - x^0\|_2^2 + h^2 M^2 k \le 2\|x_* - x^0\|_2^2,$$

т.е.

$$\|x^k - x_*\|_2 \le \sqrt{2}\|x^0 - x_*\|_2, \ k = 0,...,N. \qquad (4)$$

Для не гладких задач оценка (3) является неулучшаемой (здесь и далее неулучшаемость оценок подразумевает выполнение предположения, что размерность пространства, в котором происходит оптимизация, достаточно большая: количество итераций, которые сделает метод, не больше размерности пространства) с точностью до мультипликативного множителя [15]. Однако, если дополнительно известно, что градиент $f(\overline{x}^N)$ липшицев

$$\|\nabla f(y) - \nabla f(x)\|_2 \le L\|y - x\|_2,$$

где $x, y \in B_R(x_*)$ (см. (6)), то

$$\frac{1}{2L}\|\nabla f(x^k)\|_2^2 \le f(x^k) - f_*. \qquad (5)$$

Это неравенство является формальной записью простого геометрического факта, что если в точке $x^k$ к функции $f(x)$ провести касательную

$$f(x^k) + \langle \nabla f(x^k), x - x^k \rangle$$

и на основе такой касательной построить параболу

$$f(x^k) + \langle \nabla f(x^k), x - x^k \rangle + \frac{L}{2}\|x - x^k\|_2^2,$$

то эта парабола будет мажорировать функцию $f(x)$, т.е.

$$f(x) \le f(x^k) + \langle \nabla f(x^k), x - x^k \rangle + \frac{L}{2}\|x - x^k\|_2^2.$$

В частности выписанное неравенство имеет место и в точке минимума параболы

$$x^k - \frac{1}{L}\nabla f(x^k).$$

Поскольку приращение параболы при переходе аргумента от точки $x^k$ к точке минимума параболы составляет

$$\frac{1}{2L}\|\nabla f(x^k)\|_2^2,$$

то получаем неравенство (5). Это неравенство позволяет уточнить проведенные выше рассуждения. Как и раньше запишем

$$\|x^{k+1} - x_*\|_2^2 = \|x^k - x_*\|_2^2 - 2h\langle \nabla f(x^k), x^k - x_* \rangle + h^2\|\nabla f(x^k)\|_2^2 \le$$
$$\le \|x^k - x_*\|_2^2 - 2h(f(x^k) - f_*) + 2Lh^2(f(x^k) - f_*) = \|x^k - x_*\|_2^2 + 2h(Lh - 1)(f(x^k) - f_*).$$

При $h \le 1/L$ отсюда имеем

$$\|x^{k+1} - x_*\|_2^2 \le \|x^k - x_*\|_2^2, \ k = 0,...,N-1,$$

следовательно

$$\|x^k - x_*\|_2 \le \|x^0 - x_*\|_2, \ k = 0,...,N. \qquad (6)$$

Полагая

$$h = \frac{1}{2L},$$

получим

$$f\left(\overline{x}^N\right) - f_* \leq \frac{1}{N}\sum_{k=0}^{N-1} f\left(x^k\right) - f(x_*) \leq \frac{2L}{N}\sum_{k=0}^{N-1}\left\{\left\|x^k - x_*\right\|_2^2 - \left\|x^{k+1} - x_*\right\|_2^2\right\} \leq \frac{2LR^2}{N}.$$

**Замечание 1.** То что $f\left(\overline{x}^N\right) - f_*$ должно зависеть только от $MR/\sqrt{N}$ и(или) $LR^2/N$ является следствием П-теоремы теории размерностей [16]. Вводя точность

$$f\left(\overline{x}^N\right) - f_* \leq \varepsilon,$$

можно показать, что существует всего две (независимые) безразмерные величины, сконструированные из введенных параметров:

$$\frac{M^2R^2}{\varepsilon^2} \text{ и } \frac{LR^2}{\varepsilon}.$$

Согласно П-теореме, любая безразмерная величина должна функционально выражаться через эти две (базисные). В частности,

$$N = G\left(\frac{M^2R^2}{\varepsilon^2}, \frac{LR^2}{\varepsilon}\right).$$

В случае, когда нельзя гарантировать липшицевость градиента ситуация упрощается

$$N = \tilde{G}\left(\frac{M^2R^2}{\varepsilon^2}\right).$$

Полученная методом (2) оценка (3) соответствует этой формуле. Более того, как уже отмечалось, это оценка является неулучшаемой на классе негладких выпуклых задач. К сожалению, полученная оценка скорости сходимости метода (2) с шагом $h = 1/(2L)$ в гладком случае уже не будет оптимальной.

Кстати сказать, из П-теоремы также следует, что шаг метода (2) $h$ в негладком случае должен вычисляться по формуле

$$h = c\frac{\varepsilon}{M^2}, \ c > 0,$$

которую можно получить из использовавшейся нами ранее формулы

$$h = \frac{R}{M\sqrt{N}},$$

если выразить $N$ через $\varepsilon$ с помощью формулы (3), полагая

$$\varepsilon = \frac{MR}{\sqrt{N}}.$$

В гладком случае $h$ определяться из соотношения вида

$$W\left(h\frac{M^2}{\varepsilon}, hL\right) = 1.$$

В стохастическом случае [4] (вместо градиента получаем стохастический градиент с дисперсией $\sigma^2$) из

$$\tilde{W}\left(h\frac{M^2}{\varepsilon}, hL, h\frac{\sigma}{R}\right) = 1.$$

При условии липшицевости градиента выписанная оценка скорости сходимости может быть улучшена [15]. Например, при использовании метода сопряженных градиентов [15, 17]. Локальные оценки скорости сходимости метода тяжелого шарика также говорят об этом [17]. Однако среди большого многообразия "ускоренных методов", которые сходятся по нижним оценкам, мы особо выделим быстрый градиентный метод, предложенный в 1983 г. в кандидатской диссертации Ю.Е. Нестерова (научным руководителем был проф. Б.Т. Поляк). Помимо того, что это был один из первых методов (без использования вспомогательных одномерных и двумерных оптимизаций [15]), для которого было получено строгое доказательство глобальной сходимости согласно нижним оценкам (в гладком случае), метод оказался обладающим хорошими свойствами типа (6). Но наиболее важное свойство для нас в данной статье – его прямо-двойственность. Развитие и использование этого метода отражено в диссертации [5].

Перейдем к построению быстрого градиентного метода (БГМ). Далее мы будем во многом следовать недавно предложенному способу понимания БГМ [14]. Тем не менее, нам потребуется из этих рассуждений получить свойство типа (6) и прямо-двойственность. Наличие у БГМ этих свойств в работе [14] не исследовалось, поэтому далее мы приведем все необходимые рассуждения в нужном объеме.

Предварительно определим две числовые последовательности шагов $\{\alpha_k, \tau_k\}$:

$$\alpha_1 = \frac{1}{L}, \ \alpha_k^2 L = \alpha_{k+1}^2 L - \alpha_{k+1}, \ \tau_k = \frac{1}{\alpha_{k+1} L}.$$

Можно написать явные формулы. Нам так же пригодится упрощенный вариант этих последовательностей [14], который определяется следующим образом:

$$\alpha_1 = \frac{1}{L}, \ \alpha_k^2 L = \alpha_{k+1}^2 L - \alpha_{k+1} + \frac{1}{4L}, \ \tau_k = \frac{1}{\alpha_{k+1} L}.$$

В таком случае

$$\boxed{\alpha_{k+1} = \frac{k+2}{2L}, \ \tau_k = \frac{1}{\alpha_{k+1} L} = \frac{2}{k+2}.}$$

$$\text{БГМ}\left(x^0 = y^0 = z^0\right)$$

$$\boxed{\begin{aligned}&1. \ x^{k+1} = \tau_k z^k + (1-\tau_k) y^k; \\ &2. \ y^{k+1} = x^{k+1} - \frac{1}{L} \nabla f\left(x^{k+1}\right); \\ &3. \ z^{k+1} = z^k - \alpha_{k+1} \nabla f\left(x^{k+1}\right). \end{aligned}}$$

Из последней формулы в доказательстве леммы 4.3 [14] имеем (для всех $x$)

$$\alpha_{k+1}^2 L f\left(y^{k+1}\right) - \left(\alpha_{k+1}^2 L - \alpha_{k+1}\right) f\left(y^k\right) \le$$

$$\le \alpha_{k+1} \left\{ f\left(x^{k+1}\right) + \left\langle \nabla f\left(x^{k+1}\right), x - x^{k+1} \right\rangle \right\} + \frac{1}{2} \left\| z^k - x \right\|_2^2 - \frac{1}{2} \left\| z^{k+1} - x \right\|_2^2.$$

Просуммировав то что получается по $k = 0, \ldots, N-1$, получим

$$\alpha_N^2 L f\left(y^N\right) \le \min_x \left\{ \sum_{k=0}^{N-1} \alpha_{k+1} \left\{ f\left(x^{k+1}\right) + \left\langle \nabla f\left(x^{k+1}\right), x - x^{k+1} \right\rangle \right\} + \frac{1}{2} \left\| z^k - x \right\|_2^2 - \frac{1}{2} \left\| z^{k+1} - x \right\|_2^2 \right\} \le$$

$$\le \left( \sum_{k=0}^{N-1} \alpha_{k+1} \right) f_* + \frac{1}{2} \left\| z^0 - x_* \right\|_2^2 - \frac{1}{2} \left\| z^N - x_* \right\|_2^2. \quad (7)$$

Заметим, что

$$\max \left\{ \left\| x^k - x_* \right\|_2, \left\| y^k - x_* \right\|_2, \left\| z^k - x_* \right\|_2 \right\} \le \left\| x^0 - x_* \right\|_2, \ k = 0, \ldots, N. \quad (8)$$

Действительно, полагая в формуле (7) $N := k$, $k := i$, и учитывая, что $\alpha_i$ не зависят от $k$, а $f(y_k) \ge f_*$ и $\alpha_k^2 L = \sum_{i=0}^{k-1} \alpha_{i+1}$, получим (в точности эту же формулу можно получить и для упрощенной схемы выбора шагов)

$$\left\| z^k - x_* \right\|_2^2 \le \left\| z^0 - x_* \right\|_2^2.$$

Используя это неравенство, неравенство (6) для $y^{k+1}$ и выпуклость квадрата евклидовой нормы, получим

$$\left\| y^{k+1} - x_* \right\|_2^2 \le \left\| x^{k+1} - x_* \right\|_2^2 = \left\| \tau_k \cdot (z^k - x_*) + (1-\tau_k)(y^k - x_*) \right\|_2^2 \le$$

$$\le \tau_k \left\| z^k - x_* \right\|_2^2 + (1-\tau_k) \left\| y^k - x_* \right\|_2^2 \le \tau_k \left\| z^0 - x_* \right\|_2^2 + (1-\tau_k) \left\| y^k - x_* \right\|_2^2 =$$

$$= \tau_k \|x^0 - x_*\|_2^2 + (1 - \tau_k)\|y^k - x_*\|_2^2 = \tau_k \|y^0 - x_*\|_2^2 + (1 - \tau_k)\|y^k - x_*\|_2^2.$$

Отсюда по индукции получаем (8).

Возвращаясь к формуле (7), установим прямо-двойственность БГМ. Для этого перепишем формулу (7) для упрощенного варианта выбора шагов

$$\frac{(N+1)^2}{4L} f(y^N) + \sum_{k=0}^{N-1} \frac{1}{4L} f(y^k) \leq \min_x \left\{ \sum_{k=0}^{N-1} \frac{k+2}{2L} \left\{ f(x^{k+1}) + \langle \nabla f(x^{k+1}), x - x^{k+1} \rangle \right\} + \frac{1}{2}\|z^0 - x\|_2^2 \right\},$$

т.е.

$$f(\breve{y}^N) \leq \frac{4L}{N \cdot (N+3)} \min_x \left\{ \sum_{k=0}^{N-1} \frac{k+2}{2L} \left\{ f(x^{k+1}) + \langle \nabla f(x^{k+1}), x - x^{k+1} \rangle \right\} + \frac{1}{2}\|z^0 - x\|_2^2 \right\}, \quad (9)$$

где

$$\breve{y}^N = \frac{1}{N \cdot (N+3)} \left( \sum_{k=0}^{N-1} y^k + (N+1)^2 y^N \right).$$

Собственно именно неравенство (9) и даст нам возможность в следующем пункте восстанавливать решение прямой задачи, исходя из решения двойственной методом БГМ.

Сформулируем основной результат данного пункта.

**Теорема 1.** *Пусть функционал задачи (1) обладает свойством*
$$\|\nabla f(y) - \nabla f(x)\|_2 \leq L\|y - x\|_2, \ x, y \in B_R(x_*). \tag{10}$$

*Тогда БГМ генерирует такую последовательность точек $\{x^k, y^k, z^k\}_{k=0}^N$, что имеют место соотношения (8) и (9), причем формулу (9) можно переписать следующим образом*

$$f(\breve{y}^N) - f_* \leq \frac{2L\|z^0 - x_*\|_2^2}{N \cdot (N+3)} \leq \frac{2L\|z^0 - x_*\|_2^2}{(N+1)^2} = \frac{2LR^2}{(N+1)^2}.$$

**Замечание 2.** Новизна в этой теореме по сравнению со всеми известными ее аналогами заключается в том, что хотя задача оптимизации (1) решается на неограниченном множестве, параметры, входящие в оценки скорости сходимости, определяются расстоянием от точки старта до решения. То есть весь итерационный процесс будет находиться в (евклидовом) шаре с центром в решении и радиуса равного расстоянию от точки старта до решения. Поскольку мы априорно, как правило, не знаем это расстояние, то может показаться, что ценность в таком замечании небольшая. Однако, как будет продемонстрировано в следующем пункте, сделанное замечание играет важную роль в обосновании предлагаемого прямо-двойственного подхода (впрочем, существует другой способ рассуждений, при котором условие типа (8) не используется [18]). Интересно также заметить, что если решение задачи (1) не единственно (обозначим множество решений через $X$), то под $R^2$ можно понимать $\min_{x_* \in X} \|z^0 - x_*\|_2^2$.

**Замечание 3.** Приведенную выше теорему 1 можно распространить и на случай, когда множество, на котором происходит оптимизация, не совпадает со всем пространством. Например, является неотрицательным ортантом или симплексом. В общем случае это требует также введения прокс-структуры [5, 19] в задачу, отличной от использованной нами выше – евклидовой. Все это несколько усложняет выкладки. В частности, в описании БГМ нужно использовать шаги типа прямого градиентного метода, а их проксимальные варианты [4]. Насколько нам известно, пока это все сделано (в общности теоремы 1) только для евклидовой прокс-структуры, но с произвольными множествами, на которых ведется оптимизация.

**Замечание 4.** В действительности, описанную выше конструкцию, с сохранением основного результата – теоремы 1, можно перенести на задачи композитной оптимизации [19, 20], на задачи, в которых не известна константа $L$, и ее требуется подбирать по ходу процесса [5, 19]. Также на основе БГМ и только что сделанного замечания, можно построить и соответствующий вариант универсального метода Ю.Е. Нестерова [21]. Все описанные обобщения можно сделать и в концепции неточного оракула [4, 11, 22–25].

## 3. Приложение к задаче минимизации сильно выпуклого функционала простой структуры при аффинных ограничениях

Пусть требуется решать задачу

$$g(x) \to \min_{Ax=b,\, x \in Q}, \qquad (11)$$

где функция $g(x)$ – 1-сильно выпуклая в $p$-норме $(1 \le p \le 2)$. Построим двойственную задачу

$$F(y) = \max_{x \in Q}\{\langle y, b - Ax\rangle - g(x)\} \to \min_y. \qquad (12)$$

Во многих важных приложениях основной вклад в вычислительную сложность внутренней задачи максимизации дает умножение $Ax$ ($A^T y$). Это так, например, для сепарабельных функционалов

$$g(x) = \sum_{k=1}^{n} g_k(x_k)$$

и параллелепипедных ограничениях $Q$. В частности, это имеет место для задач энтропийно-линейного программирования (ЭЛП) [1], в которых имеется явная формула $x(y)$.

В общем случае внутренняя задача максимизации не решается точно (по явным формулам). Тем не менее, за счет сильной выпуклости $g(x)$ (аналогичное можно сказать в случае сепарабельности $g(x)$, но отсутствии сильной выпуклости) точность решения этой вспомогательной задачи (на каждой итерации внешнего метода) входит в оценку сложности ее решения логарифмическим образом, как следствие, оговорки о неточности оракула, выдающего градиент для внешней задачи минимизации $F(y)$, можно опустить. Аккуратный учет этого всего приводит лишь к логарифмическим поправкам в итоговых оценках сложности метода (см., например, [8, 19, 24, 25]). Поэтому для большей наглядности (следуя совету А.С. Немировского) мы далее в рассуждениях будем просто считать, что есть явная формула $x(y)$.

Применим к задаче (12) БГМ из п. 2 ($z^0 = 0$), получим по формуле (9) (были сделаны переобозначения: $x \to y$, $\breve{y} \to \tilde{y}$)

$$F(\tilde{y}^N) \le \frac{4L}{N \cdot (N+3)} \min_y \left\{ \sum_{k=0}^{N-1} \frac{k+2}{2L}\{F(y^{k+1}) + \langle \nabla F(y^{k+1}), y - y^{k+1}\rangle\} + \frac{1}{2}\|z^0 - y\|_2^2 \right\}, \qquad (13)$$

где (см., например, [5])

$$L = \max_{\|x\|_p \le 1} \|Ax\|_2^2.$$

В частности, для задачи ЭЛП $p = 1$ [1, 13]

$$L = \max_{k=1,\ldots,n} \|A^{\langle k \rangle}\|_2^2,$$

где $A^{\langle k \rangle}$ – $k$-й столбец матрицы $A^{\langle k \rangle}$. Для задачи PageRank (см. замечание 10 ниже) $p = 2$

$$L = \lambda_{\max}(A^T A) = \sigma_{\max}(A).$$

Из неравенства (13) имеем ($R^2 = \|z^0 - y_*\|_2^2 = \|y_*\|_2^2$, где $y_*$ – решение задачи (12))

$$F(\tilde{y}^N) - \min_{y \in B_{3R}(0)} \left\{ \sum_{k=0}^{N-1} \frac{2(k+2)}{N \cdot (N+3)}\{F(y^{k+1}) + \langle \nabla F(y^{k+1}), y - y^{k+1}\rangle\} \right\} \le \frac{18LR^2}{(N+1)^2} \stackrel{def}{=} \gamma_N. \qquad (14)$$

Введем

$$\lambda_k = \frac{2(k+2)}{N \cdot (N+3)},$$

$$x^N = \sum_{k=0}^{N-1} \lambda_k x(y^{k+1}) = \frac{2}{N \cdot (N+3)} \sum_{k=0}^{N-1} (k+2) x(y^{k+1}) = \frac{N^2+N-2}{N \cdot (N+3)} x^{N-1} + 2\frac{N+1}{N \cdot (N+3)} x(y^N).$$

Перепишем неравенство (14) исходя из определения $F(y)$ (12) и $x(y)$ (приводимая далее выкладка аналогична рассуждениям из п. 3 работы [7])

$$F(\tilde{y}^N) - \sum_{k=0}^{N-1} \lambda_k \langle y^{k+1}, b - Ax(y^{k+1}) \rangle + \sum_{k=0}^{N-1} \lambda_k g(x(y^{k+1})) -$$
$$- \min_{y \in B_{3R}(0)} \left\{ \sum_{k=0}^{N-1} \lambda_k \langle b - Ax(y^{k+1}), y - y^{k+1} \rangle \right\} \le \gamma_N.$$

Учитывая, что
$$\sum_{k=0}^{N-1} \lambda_k = 1,$$

получим
$$F(\tilde{y}^N) + g\left(\sum_{k=0}^{N-1} \lambda_k x(y^{k+1})\right) + \max_{y \in B_{3R}(0)} \left\{ \left\langle A\sum_{k=0}^{N-1} \lambda_k x(y^{k+1}) - b, y \right\rangle \right\} \le \gamma_N,$$

т.е.
$$F(\tilde{y}^N) + g(x^N) + 3R \|Ax^N - b\|_2 \le \gamma_N.$$

Отсюда (приводимые далее рассуждения во многом повторяют рассуждениями из п. 6.11 работы [3]), используя то, что
$$Ax_* = b,$$

и слабую двойственность
$$-g(x_*) \le F(y_*),$$

получаем
$$g(x^N) - g(x_*) \le g(x^N) + F(y_*) \le g(x^N) + F(\tilde{y}^N) \le g(x^N) + F(\tilde{y}^N) + 3R \|Ax^N - b\|_2 \le \gamma_N.$$

Исходя из определения $F(y)$ (12) и свойства (8), имеем
$$-g(x_*) = \langle y_*, b - Ax_* \rangle - g(x_*) = F(y_*) \ge \langle y_*, b - Ax^N \rangle - g(x^N) \implies$$
$$g(x_*) - g(x^N) \le R \|Ax^N - b\|_2,$$
$$R \|Ax^N - b\|_2 \le -g(x^N) + \langle \tilde{y}^N, b - Ax^N \rangle + g(x^N) + 3R \|Ax^N - b\|_2 \le$$
$$\le F(\tilde{y}^N) + g(x^N) + 3R \|Ax^N - b\|_2 \le \gamma_N.$$

Отсюда следует, что
$$|g(x^N) - g(x_*)| \le \gamma_N, \ R \|Ax^N - b\|_2 \le \gamma_N.$$

Поскольку
$$g(x^N) - g(x_*) \le F(\tilde{y}^N) + g(x^N) \le \gamma_N,$$

то мы приходим к следующему результату.

**Теорема 2.** *Пусть нужно решить задачу (11) посредством перехода к задаче (12), исходя из выписанных выше формул. Выбираем в качестве критерия останова БГМ выполнение следующих условий*

$$F(\tilde{y}^N) + g(x^N) \le \varepsilon, \ \|Ax^N - b\|_2 \le \tilde{\varepsilon}.$$

*Тогда БГМ гарантированно остановится, сделав не более чем*

$$\max\left\{\sqrt{\frac{18LR^2}{\varepsilon}}, \sqrt{\frac{18LR}{\tilde{\varepsilon}}}\right\}$$

*итераций.*

**Замечание 5.** Обратим внимание, что при описании подхода одновременного решения прямой и двойственной задачи мы имели неизвестный параметр $R$. Однако этот параметр не входит в описание алгоритма и в критерий его останова. Он входит только в оценку числа итераций. К сожалению, такого удается добиться далеко не всегда. Нетривиально то, что нам удалось этого добиться в данном контексте. Обычно при решении двойственной задачи искусственно компактифицируют [2, 3, 11] множество, на котором происходит оптимизация (для двойственной задачи это, как правило, либо все пространство, либо прямое произведение пространства на неотрицательный ортант). В результате используются методы, в которых требуется проектирование на шар заранее неизвестного радиуса. Эту проблему (неизвестности размера двойственного решения) обычно решают либо с помощью процедуры рестартов [1, 4, 24], что, как правило, приводит к увеличению числа итераций как минимум на порядок, либо с помощью слейтеровской релаксации, которая часто еще более затратна в смысле требуемого числа итераций [1].

**Замечание 6.** К сожалению, в ряде приложений имеет место лишь строгая (не сильная) выпуклость $g(x)$. В этом случае (хотя двойственная задача и будет гладкой) ничего нельзя сказать о константе Липшица градиента, которая явно входит в шаг БГМ. Как уже отмечалось ранее, проблема решается адаптивным подбором $L$, а в более общем случае (когда и гладкость $g(x)$ нельзя гарантировать) с использованием универсального метода [11, 21, 25]. Тем не менее, ранее открытым был вопрос, насколько все эти конструкции (описанные выше) работают при решении гладкой двойственной задачи (на неограниченном множестве) с не ограниченными равномерно константами Липшица (на этом множестве). Ведь если предположить, что метод может отдаляться от решения по ходу итерационного процесса больше, чем на расстояние, которое было в момент старта, и расстояние, на которое он может отдаляться, зависит от свойств гладкости функционала, то образуется "порочный круг". При неаккуратном оценивании так и получается. В данной статье было показано, что для детерминированных постановок задач при правильном подборе прямо-двойственных методов таких проблем можно избежать. Причем избежать естественным образом, т.е. не за счет искусственной компактификации (как это общепринято), приводящей к дополнительным затратам на рестарты.

**Замечание 7.** Класс задач, к которым применим описанный выше подход можно существенно расширить, допуская, например, в постановке задачи (11) ограничения вида неравенств $Cx \le d$, и в более общем случае вида $Cx - d \in K$, где конус $K$ имеет простое двойственное описание [26]. Можно исходить не из задач вида (11), а сразу из задачи минимизации функционала, имеющего лежандрово представления вида (12) [24]. При этом минимизация по $y$ теперь может вестись по произвольному выпуклому множеству.

**Замечание 8.** В случае, когда размерность двойственного пространства небольшая можно использовать вместо БГМ метод эллипсоидов (не требующий гладкости двойственного функционала). Этот метод также является прямо-двойственным [3]. Интересные примеры в связи с этой конструкцией возникают, когда размерность прямого пространства огромна, но при этом есть эффективный линейный минимизационный оракул, позволяющий (несмотря на огромную размерность прямого пространства) эффективно вычислять градиент двойственного функционала [9, 10, 27].

**Замечание 9.** Полезно заметить, что описанный в этом пункте подход (особенно в связи с замечанием 8) позволяет во многих важных случаях решать задачу поиска градиентного отображения, возникающую (при проектировании на допустимое множество не очень простой структуры) на каждой итерации большинства итерационных методов [3, 5, 19]. Собственно, общая идея "разделяй и властвуй" применительно к численному решению задач выпуклой оптимизации приобрела форму итерационного процесса: на каждом шаге которого решается более простая задача, чем исходная. Можно играть на том насколько сложная задача решается на каждом шаге и том, как много таких шагов надо сделать. Хороший пример здесь – это композитная оптимизация [19, 20]. Часть сложности постановки задачи (в виде композита) перенесена полностью (без линеаризации) в задачу, которую требуется решать на каждой итерации. Если композит не очень плохой, то это не сильно меняет стоимость итерации, зато может приводить к существенному сокращению необходимого числа итераций (например, для негладкого композита). Другие примеры на эту же тему имеются здесь [11, 12, 24, 25]. Общая линия рассуждений тут приблизительно такая: любое усложнение итерации, как правило, приводит к сокращению их числа, с другой стороны в стоимость одной итерации всегда входит расчет (пересчет) градиента или его (скажем, стохастического) аналога, использующегося в методе. Для детерминированных методов, когда используется градиент, основная стоимость итерации формируется как раз за счет расчета этого градиента (как правило, это умножение матрицы на вектор, т.е. порядка $\mathrm{O}(n^2)$ арифметических операций). При рассчитанном градиенте сделать шаг итерационного процесса стоит обычно не более $\mathrm{O}(n\ln(n/\varepsilon))$ арифметических операций. Таким образом, остается довольно большой зазор, в который можно занести дополнительные вычисления, перенеся еще какую-то часть сложности постановки задачи на каждый шаг, в надежде сократить число шагов. В частности, если рассматривать задачу [8, 25]

$$\frac{1}{2}\|Ax-b\|_2^2 + \mu\sum_{k=1}^{n} x_k \ln x_k \to \min_{x \in S_n(1)},$$

с достаточно большим $\mu > 0$, то необходимо учитывать сильную выпуклость энтропийного композита. Обычный БГМ в композитном варианте для этой задачи требует на каждой итерации решения некоторой почти-сепарабельной задачи. При решении такой задачи описанные в этом пункте конструкции оказываются весьма полезными. Более того, мы привели здесь эту задачу также и потому, что при определенном значении параметра $\mu > 0$ эта задача "эквивалентна" исходной задаче (11) с функционалом вида энтропии. В работе [8] обсуждаются (с содержательной стороны и со стороны практических вычислений) эти разные способы понимания одной задачи, а также приводятся некоторые конструкции (близкие к описанным в этом пункте), которые позволяют за небольшую дополнительную плату определять параметр $\mu > 0$, чтобы имело место отмеченное соответствие (см. также [25]).

**Замечание 10 (PageRank и нижние оценки).** Приведенные в теореме 2 оценки при первом взгляде могут приводить к противоречиям. Поясним это следующим примером [24].

Задача поиска такого $x^* \in \mathbb{R}^n$, что
$$Ax^* = b$$
сводится к задаче выпуклой гладкой оптимизации
$$f(x) = \|Ax-b\|_2^2 \to \min_x.$$

Нижняя оценка для скорости решения такой задачи [15] имеет вид
$$f(x^N) \geq \Omega\left(L_x R_x^2 / N^2\right), \quad L_x = \sigma_{\max}(A), \quad R_x = \|x^*\|_2.$$

Откуда следует, что только при
$$N \geq \Omega\left(\sqrt{L_x} R_x / \varepsilon\right)$$
можно гарантировать выполнение неравенства
$$f(x^N) \leq \varepsilon^2, \text{ т.е. } \|Ax^N - b\|_2 \leq \varepsilon.$$

Однако эта нижняя оценка для специальных матриц может быть улучшена. Рассмотрим задачу поиска вектора PageRank [28] ($n \sim 10^{10}$), которую мы перепишем как
$$Ax = \begin{pmatrix} P^T - I \\ 1\ldots\ldots 1 \end{pmatrix} x = \begin{pmatrix} 0 \\ 1 \end{pmatrix} = b,$$

где $I$ – единичная матрица. По теореме Фробениуса–Перрона [29] решение такой системы с неразложимой стохастической матрицей $P$ единственно и положительно $x > 0$. Сведем решение этой системы уравнений к вырожденной задаче выпуклой оптимизации
$$\frac{1}{2}\|x\|_2^2 \to \min_{Ax=b}.$$

Построим двойственную к ней задачу (поскольку система $Ax = b$ совместна, то по теореме Фредгольма не существует такого $y$, что $A^T y = 0$ и $\langle b, y\rangle > 0$, следовательно, двойственная задача имеет конечное решение)
$$\min_{Ax=b} \frac{1}{2}\|x\|_2^2 = \min_x \max_y \left\{\frac{1}{2}\|x\|_2^2 + \langle b - Ax, y\rangle\right\} = \max_y \min_x \left\{\frac{1}{2}\|x\|_2^2 + \langle b - Ax, y\rangle\right\} = \max_y \left\{\langle b, y\rangle - \frac{1}{2}\|A^T y\|_2^2\right\}.$$

Зная решение $y^*$ двойственной задачи (например, с минимальной евклидовой нормой)
$$\langle b, y\rangle - \frac{1}{2}\|A^T y\|_2^2 \to \max_y$$
можно восстановить решение прямой задачи
$$x(y) = A^T y.$$

Кроме того, если численно решать двойственную задачу БГМ, то согласно теореме 2
$$\|Ax^N - b\|_2 \leq \frac{8 L_y R_y}{N^2},$$
где
$$L_y = \sigma_{\max}(A^T) = \sigma_{\max}(A) = L_x, \quad R_y = \|y^*\|_2.$$

Кажется, что это противоречит нижней оценке
$$\|Ax^N - b\|_2 \geq \Omega\left(\sqrt{L_x} R_x / N\right).$$

Однако важно напомнить [15], что эта нижняя оценка установлена для всех $N \le n$ ($n$ – размерность вектора $x$), и она будет улучшена, в результате описанной процедуры только если дополнительно предположить, что матрица $A$ удовлетворяет следующему условию

$$L_y R_y \ll n\sqrt{L_x} R_x,$$

что сужает класс, на котором была получена нижняя оценка

$$\Omega\left(\sqrt{L_x} R_x / N\right).$$

В типичных ситуациях можно ожидать $R_y \gg R_x$, что мешает выполнению требуемого условия.

**Замечание 11.** Определенный нюансы возникают при попытке перенесения результатов данной работы на случай, когда вместо градиента доступен только стохастический градиент. Мы не будем здесь подробно на этом останавливаться (этому случаю планируется посвятить отдельную работу). Тем не менее, отметим, что для рандомизированных покомпонентных методов ответы на многие вопросы уже удалось получить [13]. В частности, задачу из замечания 10 (PageRank) можно решать прямым ускоренным покомпонентным методом или двойственным. Оценки получаются следующими (см. [13])

$$E\left[\|Ax^N - b\|_2^2\right] = E\left[f\left(x^N\right)\right] - f_* = \mathrm{O}\left(n^2 \frac{\bar{L}_x R_x^2}{N^2}\right), \quad \bar{L}_x^{1/2} = \frac{1}{n}\sum_{k=1}^{n}\|A^{\langle k \rangle}\|_2 \le 2, \quad R_x^2 \le 2, \text{ (для прямой задачи)}$$

$$E\left[\|Ax^N - b\|_2\right] = \mathrm{O}\left(n^2 \frac{\bar{L}_y R_y}{N^2}\right), \quad \bar{L}_y^{1/2} = \frac{1}{n+1}\sum_{k=1}^{n+1}\|A_k\|_2. \text{ (для двойственной задачи)}$$

При этом если матрица $P$ имеет $sn$ ненулевых элементов ($s \ll n$), то одна итерация у обоих методов в среднем требует $\mathrm{O}(s)$ арифметических операций. Без численных экспериментов (исходя из приведенных выше оценок) сложно определить, какой подход будет предпочтительнее (в основном из-за не знания $R_y$). Этот пример является, чуть ли не единственным примером, когда удается так организовать вычисления, что сполна можно использовать разреженность задачи (итерация выполняется за $\mathrm{O}(s)$). К сожалению, из-за необходимости расчета (пересчета) $x^N$, как правило, всегда приходится тратить на одну итерацию как минимум $\mathrm{O}(n)$ арифметических операций [13]. Оказывается, существует другой (более простой) способ восстановления решения прямой задачи, исходя из решения двойственной, который лучше приспособлен к возможности учета разреженности. К описанию этого способа мы сейчас и переходим.

Распространим приведенный в статье подход на случай, когда двойственный функционал в задаче (11), в свою очередь, оказывается $\mu$-сильно выпуклым (вогнутым) в 2-норме. Для этого достаточно, чтобы функционал прямой задачи имел равномерно ограниченную константу Липшица градиента, а сама прямая задача решалась бы на всем пространстве [8, 13] (т.е. нет других ограничений, кроме ограничения $Ax = b$). Такой пример нам уже встречался в замечании 10. Воспользуемся техникой рестартов (см., например, [5, 14, 25]), которая в данном случае будет иметь следующий вид (см. теорему 2; заметим, что в этой оценке мы не использовали прямо-двойственность БГМ)

$$\frac{\mu}{2}\left\|\tilde{y}^{\bar{N}} - y_*\right\|_2^2 \le F\left(\tilde{y}^{\bar{N}}\right) - F\left(y_*\right) \le \frac{2L\left\|y^0 - y_*\right\|_2^2}{\bar{N}^2}.$$

Выбирая

$$\bar{N} = \sqrt{\frac{8L}{\mu}},$$

получим, что

$$\left\|\tilde{y}^{\bar{N}} - y_*\right\|_2^2 \le \frac{1}{2}\left\|y^0 - y_*\right\|_2^2.$$

Выберем в БГМ в качестве точке старта $\tilde{y}^{\bar{N}}$, и снова сделаем $\bar{N}$ итераций, и т.д. Несложно понять, что если мы хотим достичь точности по функции $\varepsilon$, то число таких рестартов (перезапусков) БГМ достаточно взять (здесь используется достаточное стандартное обозначение $\lceil \cdot \rceil$, которое мы поясним примером $\lceil 0.2 \rceil = 1$)

$$\left\lceil \log_2\left(\frac{\mu R^2}{\varepsilon}\right) \right\rceil.$$

Таким образом, общее число итераций, которое сделает БГМ, можно оценить, как

$$\sqrt{\frac{8L}{\mu}} \left\lceil \log_2\left(\frac{\mu R^2}{\varepsilon}\right) \right\rceil.$$

На первый взгляд, кажется, что при таком подходе мы не контролируем

$$\|Ax(y) - b\|_2 = \|\nabla F(y)\|_2.$$

В действительности, мы всегда (не только для прямо-двойственных методов) для задач с Липшицевым градиентом можем контролировать $\|\nabla F(y)\|_2$, используя неравенство (5) (если $F(y_*) \neq 0$, то вместо градиента $\nabla F(y)$ можно вычислять градиентное отображение, см., например, [4, 5])

$$\frac{1}{2L}\|\nabla F(y)\|_2^2 \leq F(y) - F(y_*).$$

Проблема, однако, в том, что типично это неравенство довольно грубое. Действительно, в $\mu$-сильно выпуклом случае

$$\frac{1}{2L}\|\nabla F(y)\|_2^2 \leq F(y) - F(y_*) \leq \frac{1}{2\mu}\|\nabla F(y)\|_2^2.$$

Если $L/\mu \gg 1$ (что типично), то использование неравенства (5) может приводить (и приводит [13]) к сильно завышенным оценкам. Однако, если наряду с использованием оценки (5), учитывать геометрическую скорость сходимости БГМ в виду наличия сильной выпуклости, то (считаем, что $\mu R^2 \geq \varepsilon$)

$$\|Ax(y^N) - b\|_2 = \|\nabla F(y^N)\|_2 \leq \sqrt{2L \cdot (F(y^N) - F(y_*))} \leq \sqrt{2L\mu R^2} \exp\left(-\frac{N}{2}\sqrt{\frac{\mu}{8L}}\right),$$

где $\{y^N\}$ обозначает последовательность, которую генерирует описанный выше БГМ с рестартами.

Перейдем к описанию критерия останова метода БГМ с рестартами. По определению $x(y)$ имеем

$$g(x(y)) + \langle y, Ax(y) - b \rangle \leq g(x_*).$$

Откуда

$$g(x(y)) - g(x_*) \leq \|y\|_2 \|Ax(y) - b\|_2.$$

Таким образом, критерий останова будет иметь вид

$$\|y^N\|_2 \|Ax(y^N) - b\|_2 \leq \varepsilon, \quad \|Ax(y^N) - b\|_2 \leq \tilde{\varepsilon}. \tag{15}$$

**Теорема 3.** *Пусть нужно решить задачу (11) посредством перехода к задаче (12) с $\mu$-сильно выпуклым в 2-норме функционалом, исходя из выписанных выше формул. Выбираем в качестве критерия останова (15). Тогда описанный выше БГМ с рестартами гарантированно остановится, сделав не более чем*

$$\max\left\{\sqrt{\frac{8L}{\mu}}\left\lceil\log_2\left(\frac{2L\mu R^4}{\varepsilon^2}\right)\right\rceil, \sqrt{\frac{8L}{\mu}}\left\lceil\log_2\left(\frac{2L\mu R^2}{\tilde{\varepsilon}^2}\right)\right\rceil\right\}$$

*итераций.*

**Замечание 12.** В действительности, только что мы описали довольно общий способ решения большого числа задач, путем перехода к гладкой двойственной. Все что нужно добавить к описанному выше подходу – это искусственное введение регуляризации $\mu \approx \varepsilon/R^2$ в двойственную задачу, когда она не сильно выпуклая [25]. При этом возникают рестарты по параметру $\mu$, поскольку $R$ априорно не известно [1, 4]. Собственно, в [1] именно таким образом и предлагалось решать задачу ЭЛП. Теоретические оценки говорят, что оба метода должны работать приблизительно в одно время (метод из настоящей работы чуть по лучше), но эксперименты [18] совершенно однозначно показали, что метод, который изложен в этой статье, работает на несколько порядков быстрее, чем метод из [1]. Связано это с тем, что из-за рестартов по $\mu$, мы обязаны сделать (на каждом шаге рестарта) предписанное число итераций метода (меньше нельзя), в то время как предложенный в данной работе метод, во-первых, не требует рестартов (это сразу экономит почти один порядок), а во-вторых, останавливается по более гибкому критерию (см. теорему 2), допускающему, что может быть сделано меньше итераций, чем предписано полученной оценкой. Численные эксперименты показывают, что именно в этом месте и происходит основная экономия нескольких порядков в оценке общего времени работы нового метода.

**Замечание 13.** Все, что описано выше для сильно выпуклого случая, переносится на произвольные методы (прямо-двойственность при этом не нужна). Например, на метод сопряженных градиентов или метод Ньютона и их модификации [17] (есть общий тезис принадлежащий А.С. Немировскому, что практически любой разумный численный метод либо уже является прямо-двойственным, либо имеет соответствующую модификацию – для многих важных методов их прямо-двойственность уже установлена в различных работах, однако, насколько нам известно, для метода сопряженных градиентов и метода Ньютона этого пока не было сделано). В частности, критерий останова (15) является общим способом контроля точности решения, полученным произвольным методом, с помощью которого хотим одновременно решить прямую и двойственную задачу. Как уже отмечалось выше в случае, когда двойственная задача только гладкая (это необходимо для справедливости проведенных рассуждений), но не сильно выпуклая, то имеющиеся сейчас теоретические способы оценивания скорости сходимости с критерием останова (15) типично приводят к сильно завышенным оценкам. Это не означает, что соответствующие методы будут плохими, просто сейчас, насколько нам известно, не существует точных способов оценивания. Проблема наличия теоретического обоснования решается путем регуляризации двойственной задачи (см. замечание 12).

**Замечание 14.** Описанный выше подход получения БГМ для сильно выпуклых задач, базирующийся на рестартах по расстоянию от текущей точки до решения, имеет один существенный недостаток. На каждом рестарте необходимо сделать предписанное число итераций, которое, как правило, оказывается завышенным. Выйти из цикла раньше возможно, если есть критерий останова. Но $F(y_*)$, как правило, не известно. На практике можно попробовать контролировать малость нормы градиента (в общем случае нормы градиентного отображения), и исходить из уменьшения квадрата этой нормы в два раза, однако для такого подхода пока не доказано, что сохраняются аналогичные (неулучшаемые) по порядку оценки итоговой скорости сходимости. Из этой ситуации есть выход – использовать вместо БГМ с рестартами БГМ без рестартов для сильно выпуклых задач [5]. А еще лучше БГМ без рестартов для сильно выпуклых задач с адаптивным подбором константы Липшица градиента. Такой метод (к тому же непрерывный по параметру сильной выпуклости) описан, например, в работе [30]. Если параметр сильной выпуклости неизвестен (см. замечание 12), то, к сожалению, все равно не удается уйти от рестартов, но рестарты уже будут только по этому параметру, и выход из последнего рестарта (самого дорого, вносящего основной вклад в оценку общего времени работы метода) может быть осуществлен (в отличие от описанного выше подхода) по контролю малости нормы градиента (градиентного отображения). Другой способ усовершенствования конструкции рестартов изложен в работе [31].